\documentclass[11pt]{article}
\def\ti{\tilde}
\def\f{\frac}
\def\0{\bf 0}
\def\1{\bf 1}
\def\t{\theta}
\def\l{\lambda}
\def\b{\beta}
\def\bd{\buildrel}
\def\Pb{{\bf P}}
\def\s{\sigma}
\def\Nb{\mathbf{N}}
\def\Lb{{\bf L}}
\def\Rb{{\bf R}}
\def\Cb{{\bf C}}
\def\Cov{{\bf Cov}}
\def\Eb{{\bf E}}
\def\E{{\bf E}}
\def\g{\gamma}
\def\iy{\infty}
\def\vn{\vskip 0.5cm\n}
\def\nl{\par\n}
\def\n{\noindent}
\def\inty{\int_{-\iy}^\iy}
\def\hs{\hskip 0.5cm}
\def\rt{\right}
\def\lt{\left}
\def\e{\varepsilon}
\def\Ec{{\cal E}}
\def\h{\hat}
\def\Fc{{\cal F}}
\def\Mc{{\cal M}}
\def\sun{\sum_{i=1}^n}
\def\endsymbol{$\sqcup\mkern-12mu\sqcap$}
\def\done{\ \endsymbol\medskip}
\def\zetb{\mbox{\boldmath$\zeta$}}
\def\etab{\mbox{\boldmath$\eta$}}
\def\Var{\mathop{\rm\bf Var}\nolimits}
\usepackage{amsmath,amsfonts}

\begin{document}

\title{On estimation of analytic density function in $L_p$}
\author{ Natalia Stepanova\\ Carleton University,
Canada}

\maketitle

\begin{abstract}
Let $X_1,X_2,\ldots$ be a sequence of independent identically distributed random variables with an unknown density function $f$ on $\Rb$.
The function $f$ is assumed to belong to a certain class of analytic functions.
The problem of estimation of $f$ using $\Lb_p$-risk, $1\leq p<\iy$, is considered.
A kernel-type estimator $f_n$ based on $X_1,\ldots,X_n$ is proposed and the upper bound on its limiting local minimax risk is established.
Our result is consistent with a conjecture of
 Guerre and Tsybakov (1998) and augments previous work in this area.

\end{abstract}

 \section{Introduction} Let $X_1,X_2,\ldots$ be a sequence of
i.i.d. random variables with unknown density function $f_0$ on $\Rb$. The function $f_0$
is assumed to belong to a  class of analytic functions. The problem
is to estimate $f_0$ using $\Lb_p$-risk, $1\leq p<\iy$.
Specifically, to judge the quality of an estimator
${f}_n(x)=f_n(x,X_1,\ldots,X_n)$ of $f_0(x)$ we will be using  $\Lb_p$-risk, $1\leq p<\iy$,
based on a loss function $L(x)$ of the type
\begin{gather*}
L(x)=l(\|x\|_p),
\end{gather*}
where $l:[0,\iy)\to \mathbb{R}$ is a function from a general class of loss functions ${\cal L}$ and for a function $f$ in $\Lb_p=\Lb_p(\Rb)$
$$\|f\|_p=\left(\int_\Rb|f(t)|^p\, dt\right)^{1/p}, \quad 1\leq p<\iy.$$

In nonparametric estimation, we only know a prior that an unknown function (such as density function, distribution function, or regression function)
belongs to some functional class. Let $f$ be an estimated function that is known to belong to a functional space ${\cal F}$.
For $l\in{\cal L}$  consider the  minimax $\Lb_p$-risk
\begin{gather*}
{\cal R}_n({\cal F})=\inf_{f_n}\sup_{f\in{\cal F}}\E_f l(\|f_n-f\|_p),
\end{gather*}
where $f_n$ is an arbitrary estimator of $f$.
If an estimator  $\ti{f}_n$ of  $f$ is such that
\begin{gather*}
\sup_{f\in{\cal F}}\E_f l(\|\ti{f}_n-f\|_p)\sim {\cal R}_n({\cal F}),
\end{gather*}
where the relation $a_n\sim b_n$ means $\lim_{n\to \iy}a_n/b_n=1$,
then $\ti{f}_n$ is called an  \textit{asymptotically efficient} (or \textit{asymptotically minimax}) estimator of $f$ with respect to $\Lb_p$-risk.
After the work of Pinsker (1980) the problem of constructing asymptotically efficient  estimators with known \textit{exact rates},
 i.e., optimal rates including optimal constants, aroused considerable interest.
 In this regard, we can mention the works of Nussbaum (1983), Ibragimov and Hasminskii (1984), Korostelev (1993),
 Donoho (1994), Golubev, Levit, and Tsybakov (1996), Schipper (1996), Guerre and Tsybakov (1998), Levit and Stepanova (2004), etc.

Sometimes estimated functions admit more precise locally asymptotically efficient estimators (see, e.g., Golubev and Levit (1996), Belitser (1998)).
Specializing to the estimation of the density function $f_0$,
suppose that for any sufficiently small vicinity ${\cal V}$ of $f_0$
in an appropriate topology, an estimator  $\ti{f}_n$ satisfies
\begin{gather}\label{krit}
\sup_{f\in{\cal V}}\E_f l(\|\ti{f}_n-f\|_p)\sim \inf_{f_n}\sup_{f\in{\cal V}}\E_f l(\|f_n-f\|_p).
 \end{gather}
 Then $\ti{f}_n$ is called a  \textit{locally asymptotically efficient} (or \textit{locally asymptotically minimax}) estimator of $f_0$ with respect to $\Lb_p$-risk.

  Among a variety of functional classes traditionally studied in nonparametric estimation, a class of analytic functions plays an important role (see, e.g., Levit and Stepanova (2004), pp.\,254--255, for discussion and references). As an alternative to classes of functions of finite smoothness,
  analytic functions  were first used for the purpose of nonparametric estimation by Ibragimov and Hasminskii (1980).
For the famous Gaussian white noise model, the problem of asymptotically efficient estimation
of analytic regression function using $\Lb_p$-risk, $1\leq p<\iy$, was solved
by Guerre and Tsybakov (1998) in a univariate case, and by Levit and Stepanova (2004) in a multivariate case.
At the same time, due to the lack of normality of observations,  an analogous problem of constructing (locally) asymptotically efficient estimators  of
analytic density remains largely unsolved.
A conjecture in Remark 5 of
 Guerre and Tsybakov (1998),
which is based on the asymptotics of the stochastic part of the $\Lb_p$-error of kernel density estimators in Cs\"{o}rg\H{o}
and Horv\'{a}th (1988), says that for a general class of loss functions $l$ and any $2\leq p<\iy$
\begin{gather}\label{GT}
\lim\limits_{n\to\iy}\inf_{\ti{f}_n}\sup_{f\in\mathbb{A}}\Eb_f l
(\psi_p^{-1}(n)\|\ti{f}_n-f\|_p)= \lim\limits_{n\to\iy}\sup_{f\in\mathbb{A}}\Eb_f l
(\psi_p^{-1}(n)\|{f}^*_n-f\|_p)=l(1),
\end{gather}
where $f^*_n(x)$ is the kernel density estimator of the form
\begin{gather}\label{yadro}
f^*_n(x)=\frac{1}{n h_n} \sum_{i=1}^n K\left(\frac{x-X_i}{h_n}\right),\quad K(t)=\frac{\sin t}{\pi t},\quad h_n=\frac{2\g}{\log n},
\end{gather}
$\psi_p(n)$ is the rate function given by (\ref{ratef}), and $\mathbb{A}={\mathbb{A}}(\g,M)$ is a class of density functions such that each $f\in\mathbb{A}$ admits an analytic continuation to the strip
$S_\g=\{x+iy:|y|\leq \g\}$ with $\g>0$ such that $f(x+iy)$ is analytic on the interior of $S_\g$, bounded on $S_\g$ and for some $M>0$
$$\int_{\Rb}|f(x+i\g)|^2 dx\leq M.$$
The only case when the limiting relation (\ref{GT}) is proved is for $p=2$ (see Theorem 2 of Schipper (1996)).
A partial solution to the conjecture of Guerre and Tsybakov (1998) was recently provided by Mason (2009) who showed that
for a general class of loss functions $l$, for each $2\leq p<\iy$ and $f\in\mathbb{A}(\g,M)$ one has (see Mason (2009), Sec.\,3.3)
\begin{gather*}
\lim_{n\to \iy}\E_f \left(\psi_p^{-1}(n)\|f^*_n-f\|_p\right)=l(1),
\end{gather*}
where $\psi_p(n)$ is given by (\ref{ratef}) and $f^*_n(x)$ is as in (\ref{yadro}).
We also know that for any  $2\leq p<\iy$ (see Mason (2009), p.\,102)
\begin{gather}\label{mason}
\limsup_{n\to \iy} \sup_{f\in \Lb_{p/2}(B)\cap\mathbb{A}(\g,M)}\E_f \left(\sqrt{n h_n}\|f^*_n-f\|_p\right)<\iy,
\end{gather}
where $\Lb_{p/2}(B)$ is the class of densities $f$ on $\Rb$ such that $\|f\|_{p/2}\leq B$ for some $B>1$.

The main result of this paper, Theorem 1, augments previous work on this topic. For a general class of loss functions $l$, we construct a kernel-type estimator $f_n$ of analytic density $f_0$ such that for any $1\leq p<\iy$, cf. (\ref{GT}) and (\ref{mason}),
\begin{gather}\label{uubb}
\lim_{{\cal V}\searrow f_0}\limsup\limits_{n\to\iy}\sup_{f\in\cal{V}}\Eb_f l
(\psi_p^{-1}(n)\|{f}_n-f\|_p)\leq l(1),
\end{gather}
where $\psi_p(n)$ is given by (\ref{ratef}) and ${\cal V}$ is a vicinity of  $f_0$ in a suitable topology.
 We expect that the lower bound
\begin{gather*}
\lim_{{\cal V}\searrow f_0}\liminf\limits_{n\to\iy}\inf_{\ti{f}_n}\sup_{f\in\cal{V}}\Eb_f l
(\psi_p^{-1}(n)\|\ti{f}_n-f\|_p)\geq l(1),
\end{gather*}
which jointly with (\ref{uubb}) would ensure that our estimator $f_n$ is locally asymptotically efficient,
also holds. The derivation of the lower bound, however, is a more intricate problem that
  ``requires a considerable amount of ingenuity and special techniques" (see Mason (2009), p.\,70). We shall not treat this problem here.

The paper is organized as follows. In Section 2 we define the class of analytic functions. In Section 3
we introduce the kernel-type estimator of an unknown analytic density function and establish its basic properties.
 The main result of the paper, Theorem 1, is stated in Section 4, and its proof is given in Section 5.
 Auxiliary results and their proofs are collected in Sections 6 and 7.

\section{An analytic functional class}
\subsection{Preliminaries} In order to construct an
appropriate functional class, we start with the class $\mathcal{F}(\g)$ of
functions $f(z),\ z=x+iy\in\Cb,$\ such that\vn a) $f$ is analytic in
the strip $S_\g=\{x+iy: |y|<\g\}$; \nl b) $u(x,y)=u_f(x,y)={\rm Re} f(x+iy)$
is bounded in $S_\g$; \nl c) $f$ is real on the real axis. \vn Under
these conditions the limits
$$u(x)=\lim_{y\to \g} u(x,y)=\lim_{y\to -\g} u(x,y)$$
are known to exist and be equal, for almost all $x$, and the
function $f(z)$ admits the following representation
\begin{gather}\label{eqn1}
f(z)=\inty G(z-s)u(s)\,ds,\hs |{\rm Im} z|<\g,
\end{gather} where
$$G(z)=\f 1{2\g\cosh\f{\pi z}{2\g}}.$$
Moreover, the whole class $\mathcal{F}(\g)$ will be obtained in this way by
allowing any function $u$ with $ \|u\|_\iy<\iy$ in (\ref{eqn1}). See Akhiezer
(1965), Sec.\,110, p.\,267--268 (without detailed proof); Timan
(1960), Sec.\,3.8.5, p.\,150. 

More general  analytic classes are obtained by using the following
arguments. Note that for $|y|<\g$,
$${\rm Re}\,G(x+iy)=\f1{2\g}\f {\cosh \f{\pi x}{2\g}\cos \f{\pi y}{2\g}}{\sinh^2\f{\pi x}{2\g}+\cos^2\f{\pi y}{2\g}}>0,$$
and
$${\rm Im}\,G(x+iy)=-\f1{2\g}\f {\sinh \f{\pi x}{2\g}\sin \f{\pi y}{2\g}}{\sinh^2\f{\pi x}{2\g}+\cos^2\f{\pi y}{2\g}}.$$
Thus by denoting
$$a=\cos\f{\pi y}{2\g}\hs{\rm and}\hs t=\sinh\f{\pi x}{2\g},
$$
we get for any $|y|<\g$
$$\lt\| {\rm Re}\,G(\cdot+iy)\rt\|_1=\inty {\rm Re}\, G(x+iy)\,dx=\inty\f{dt}{a\pi\lt(( t/a)^2+1\rt)}=1
$$
or equivalently
$$\lt\| {\rm{Re}}\, G(\cdot+iy)\rt\|_1=1.
$$
Therefore, by the generalized Minkowski inequality (Akhiezer (1965),
Sec.\,5, p.\,14), for any $1\le p<\iy$,
\begin{gather*}
\|{\rm Re}\, f(\cdot+iy)\|_p=\|u(\cdot,y)\|_p\le \|u\|_p,\hs
(|y|<\g).
\end{gather*}
 Assuming $\|u\|_p<\iy$, $1\leq p<\iy$, we get for $y\to\pm\g,$
cf. Stein (1970), Theorem 2b-c), pp.\,62--63,
$$u(x,y)\hs \to\hs u(x),\hs\rm a.e.
$$
and
$$\|u(\cdot,y)-u(\cdot)\|_p\hs\to\hs 0.
$$

\n These considerations lead to the following definition. For
given $\g>0$, $M>0$, and $1\le p<\iy$, denote by $\mathcal{F}(\g,p,M)$
the class of  functions $f(x),\ x\in \Rb,$ satisfying the following
conditions\nl \vn a)   $f$ admits representation (\ref{eqn1}), in
which\nl b) $\|u\|_p\le M.$ \vn Note that such functions are real on
the real line,  admit analytic continuation in the strip
$S_\g=\{z:|{\rm{Im}}\,z|<\g\}$, and for almost all $x$ have limiting
values $u(x)$ at the boundary $|{\rm Im}\,z|=\g$.

\vn \textbf{Remark 1.}
\textit{In the case of $2\pi$-periodic analytic functions $f(z)$, which are
real on the real axis, the representation (\ref{eqn1}) holds if and
only if (see Sarason (1965), Wilderotter (1996))
$$\sup_{|y|<\g}\lt(\int_0^{2\pi}|u(x,y)|^p\,dx\rt)^{1/p}<\iy.$$
It seems plausible that in our case a similar role is played by the
condition}
$$\sup_{|y|<\g}\|u(\cdot,y)\|_p<\iy.
$$

\subsection{The best harmonic approximation of classes $\mathcal{F}(\g,p,M),$
$1\leq p\leq \iy.$} Recall that an entire function $g(z)$, $z\in\Cb$, is said to be
the \textit{entire function of exponential type $N$} if for any $\e>0$ there is a positive constant $A=A_\e$
such that $|g(z)|\leq Ae^{(N+\e)|z|}$ for all $z\in\Cb$.
Denote by $\Ec_{N}$ the class of entire functions
$g(x)$ of exponential type $N$ such that $g\in\Lb_1.$ It is
known that such functions are necessarily bounded on $\Rb$ (see
Akhiezer (1965), Sec.\,83). Hence, they belong to $\Lb_p$, for any $p\ge 1.$
The Fourier transform of $g\in\Ec_N$ will be denoted by $\h g$:
$$\hat{g}(t)=\int_{\Rb} g(x) e^{-i t x} dx,\quad t\in\Rb.$$
 By the
Wiener-Paley Theorem (see Akhiezer (1965), Sec.\,82) the Fourier transform
$\h g(t)$ of any function $g\in \Ec_N$ vanishes outside the interval
$[-N,N].$

The problem of the best approximation in $\Lb_p$ of functions
$f\in \mathcal{F}(\g,p,M)$ by functions $g\in\Ec_N$ is discussed in
Akhiezer (1965), Ch.\,IV--V. Further references can be found therein. The
importance of such approximations in nonparametric density
estimation has been demonstrated in Ibragimov and Hasminskii
\,(1980).

A starting point in finding such an approximation is to find a best
harmonic approximation $G_N\in \Ec_N$ in $\Lb_1$ for the function
$G(x)$ appearing in (\ref{eqn1}). Here we cite the following result
from Akhiezer (1965), Sec.\, 88, p.\,207, and Sec.\,110, p.\,268. There
exists a function $G_N\in\Ec_N$ such that
$$\|G-G_N\|_1=\inf_{g\in \Ec_N}\|G-g\|_1=\f4\pi\sum_{k=0}^\iy \f {(-1)^k}{(2k+1)\cosh (2k+1)N\g}<\f 8\pi e^{-\g N}.
$$
Note that since $\|G_N\|_1<\iy$, the function $f^*_N(x)=(G_N*u)(x)$
is in $\Lb_p.$ Moreover,
 Akhiezer (1965) shows (see \,Sec.\,97, p.\,228; and Sec.\,100, p.\,237) that  $f^*_N\in\Ec_N.$ Thus another application of the generalized
Minkowski inequality gives, uniformly over $\mathcal{F}(\g,p,M)$,

$$\|f-f^*_N\|_p=\|(G-G_N)*u\|_p\le\|G-G_N\|_1\|u\|_p\le \f {8M}\pi e^{-\g N}.
$$
\vn This is of course a merely existence result, as $f^*_N$ cannot be
viewed as `constructive' approximation.  Indeed, finding $f^*_N$
involves knowing $u$, which in turn requires solving the integral
equation (\ref{eqn1}). This is known as {\it incorrect problem} in
the Theory of Integral Equations. Nevertheless, the very existence
of $f^*_N$ with the above properties can be used to produce a
`constructive' approximation to function $f\in\Fc(\g,p,M)$ which is
nearly as good as its best harmonic approximation. 

\subsection{Multipliers} Consider the class $\Mc_N$ of functions $g$
such that\vn a) $g\in\Lb_1(\Rb)$;\nl b) $\h g(t)\equiv 1,$ for
$|t|\le N$.\vn Such functions we will call {\it multipliers}, in a
slight variation of how this term is defined in the Theory of
Approximation. Multipliers produce {\it constructive approximation}
to a given function $f$ in any $\Lb_p, 1\le p \leq\iy$, in the form
$g*f$, which is almost as good as the best harmonic approximation of
$f$. Indeed, using the best harmonic approximation $f^*_N$ from the
above, and noting that by the Convolution Theorem $g*f^*_N\equiv
f^*_N$, we obtain, by another application of the generalized
Minkowski inequality, that for any $g\in\Mc_N$, uniformly over $\mathcal{F}(\g,p,M)$,
\begin{eqnarray}\label{eqn2}
\|g*f-f\|_p&=&\|g*(f-f_N^*+f_N^*)-f\|_p=\|g*(f-f_N^*)+(f_N^*-f)\|_p
\nonumber\\&\le&
\|g*(f-f_N^*)\|_p+\|f_N^*-f\|_p\le(1+\|g\|_1)\|f_N^*-f\|_p\nonumber\\
&\le& \f {8M(1+\|g\|_1)}\pi \,e^{-\g N}.
\end{eqnarray}
In other words, an application of a multiplier to a function $f$
affects the rate of the best harmonic approximation only by a factor
$(1+\|g\|_1).$

 \subsection{A class of density functions}
Now we introduce the functional class of interest.
When estimating density function using $\Lb_p$-risk, we  will distinguish between the two cases:
$1\leq p<2$ and $2\leq p<\iy.$ The case of $1\leq p<2$  requires separate consideration and
more severe assumptions on the underlying class of functions.
Specifically, let $\g>0$ and $M>0$ be given numbers, and let for $2\leq p<\iy$

$$\mathbb{F}(\g,p,M)=\{f\,:\, f\;\mbox{is a density on}\; \mathbf{R}\;\mbox{and} \;f\in\mathcal{F}(\g,p,M)\cap \Lb_{p/2}\},$$

\medskip
\n and for $1\leq p<2$ and some $(2-p)/p<\l\leq 2$
$$\mathbb{F}(\g,p,M)=\{f\,:\, f\;\mbox{is a density on}\; \mathbf{R}, \;f\in\mathcal{F}(\g,p,M)\;\mbox{and}\;\int_{\Rb}|x|^{\l} f(x) dx<\iy\}.
$$

\medskip

\n By Lemma 1 of Mason (2009), when $1\leq p<2$ any function $f\in\mathbb{F}(\g,p,M)$ is in $\Lb_{p/2}$.
Also, for all $1\leq p<\iy$ any function $f\in\mathbb{F}(\g,p,M)$  is a bounded function that belongs to $\Lb_p$.
Indeed, by (\ref{eqn1}) and H\"{o}lder's inequality  for any  $f\in\mathbb{F}(\g,p,M)$ and all $x\in\Rb$
\begin{eqnarray*}
|f(x)|&\leq& \int_{\Rb}|G(x-s) u(s)|ds\leq \|G\|_q \|u\|_p \\ &\leq&
\frac{M}{(2\g)^{1/p}\,\pi^{1/q}}\left(\int_0^{\iy} \frac{dx}{\cosh^q
x}\right)^{1/q}=:C(\g,p,M)<\iy,
\end{eqnarray*}
and
$$\int_{\Rb}|f(x)|^p dx\leq\max_{x\in\Rb}|f(x)|^{p/2}\int_{\Rb}|f(x)|^{p/2} dx<\iy.$$

\medskip
\n Now let $\mathbb{T}(\g,p,M)$
be the topology on $\mathbb{F}(\g,p,M)$ induced by the distance

$$\rho_p(f,g)=\left\{\begin{array}{ll} \|u_f-u_g\|_p+\int_{\Rb}|x|^{\l}|f(x)-g(x)|dx,\quad 1\leq p<2,\\ \vspace{-0.3cm}\\
\|u_f-u_g\|_p+\|f-g\|_{p/2},\quad\quad \quad\quad \quad\quad 2\leq p<\iy,
\end{array}\right.$$

\vspace{0.2cm}
\n where $\l$ is the same as in the definition of $\mathbb{F}(\g,p,M)$ and
$$u_f(x)=\lim_{y\to \pm\g}{\rm Re}\, f(x+iy),\quad u_g(x)=\lim_{y\to \pm\g}{\rm Re}\, g(x+iy).$$
This is a strong topology in the sense that closeness with respect to $\rho_p$ implies, via (\ref{eqn1}) and H\"{o}lder's inequality,
closeness in the uniform topology. Also, by the inverse triangular inequality, closeness with respect to $\rho_p$ implies closeness of the $\Lb_{p/2}$-norms, $2\leq p<\iy$.
With respect to the topology $\mathbb{T}(\g,p,M),$ we have
 \begin{enumerate}
\item[{\bf (A1)}] for all $1\leq p<2$ and some $(2-p)/p<\l\leq2$ as above, locally uniformly in $f\in  \mathbb{F}(\g,p,M)$,
the integral
 $\int_{|x|>B}|x|^{\l}f(x)dx$
converges to zero as $B\to \iy$;

\item[{\bf(A2)}] for $2\leq p<\iy$, locally uniformly in $f\in  \mathbb{F}(\g,p,M)$, the integral
$\int_{|x|>B}f^{p/2}(x)dx$
converges to zero as $B\to \iy$ .
\end{enumerate}

\medskip
\n \textbf{Remark 2.} \textit{For $1\leq p<2$, due to {\rm \textbf{(A1)}} and Lemma 1 of Mason (2009), locally uniformly in   $f\in  \mathbb{F}(\g,p,M)$,
the integral $\int_{|x|>B}f^{p/2}(x)dx$ converges to zero as $B\to \iy$.}
\bigskip

\section{Kernel-type estimators}
In this section, we assume that $X,X_1,X_2,\ldots$ are independent random variables with common density function $f\in\mathbb{F}(\g,p,M)$, $1\leq p<\iy$.
\subsection{Construction of the kernel}
To estimate analytic density with $\Lb_p$-risk, we will be using a kernel-type estimator of the form
$$ f_n(x)=\f1n\sun g(x-X_i),\quad x\in\Rb,$$
where $g\in\Mc_N$ (see Section 3.1 for the definition of ${\Mc}_N$). We will see that the systematic error (the bias
term) of a such estimator is completely  determined by $\|
g*f-f\|_p$, whereas its stochastic term is determined by
$$\|g\|_2\|f\|_{p/2}.$$
Note that by the Parseval inequality for any $g\in\Mc_N$,
$$\|g\|_2^2=(2\pi)^{-1}{\|\h g\|_2^2}\ge {N}/{\pi}.$$
Thus our goal is to find a multiplier $g\in\Mc_N$ whose
$\Lb_2$-norm is `close' to the lowest possible bound $\sqrt{N/\pi}$
and which, at the same time, has a reasonably well-behaved
$\Lb_1$-norm $\| g\|_1$. In other words, the question is --
informally -- how to make $\|g\|_2$ and $\|g\|_1$ both small? For
this,  we will consider a classical family of kernels well known in
Approximation Theory:
$$k(x;\t)=\f{\cos \t x-\cos x}{\pi(1-\t)x^2},\hs 0\le \t< 1.$$
In Akhiezer (1965) they are called {\it Fej\'er-type kernels.} For
$\t=1$ we obtain, as a limiting case, the sinc function:
$$k(x)=k(x;1)=\f{\sin x}{\pi x}.$$
As other special cases we obtain the Fej\'er $(\t=0)$ and the Vallee
Poussin ($\t=1/2$) kernels (see the Table below).
\vn {\centerline{\begin{tabular}{|c|c|c|}\hline ${
k(x)}$&${\mbox{$\theta$}}$&{ Kernel}\\  \hline
$\frac{2\sin^2\f x2}{\pi x^2}$&$0$&Fej\'er\\\hline $\frac{2(\cos(
x/2)-\cos x)}{\pi x^2}$&$1/2$&Vall\'ee-Poussin\\\hline $\frac {\sin
x}{\pi x}$&$1$&sinc function\\\hline
\end{tabular}}}
\bigskip

\n The Fourier transforms of these kernels are given by
\begin{gather}\label{ftr}
\h k(t;\t)=\mathbb{I}(|t|\le \t)+\f{1-|t|}{1-\t}\mathbb{I}(\t\le
|t|\le 1).
\end{gather}
We have
$$\|k(\cdot;\t)\|_2^2=\f{1+\t}{2\pi}$$
and, cf. Akhiezer (1965), Sec.\,106, p.\,255,
$$\f 4{\pi^2}\log\f{1+\t}{1-\t}+\f13\le\|k(\cdot;\t)\|_1\le \f 4{\pi^2}\log\f{1+\t}{1-\t}+2.$$
\vn Using the above family of kernels, it is easy to construct, for
any $0<\t<1$, a multiplier $k_N\in \Mc_N$ as
$$k_N(x;\t)=\f N\t k\lt(\f{Nx}\t;\t\rt).$$
Note that by the above formulas
$$\|k_N(\cdot;\t)\|_1=\|k(\cdot;\t)\|_1\le \f 4{\pi^2}\log\f{1+\t}{1-\t}+2$$
and
\begin{gather}\label{nomerr}
\f N\pi\le\|k_N(\cdot;\t)\|_2^2=\f N\t\|k(\cdot;\t)\|_2^2=\f N{2\pi}\lt(1+\f 1\t\rt).
\end{gather}
\vn Now, a meaningful choice of the kernel from the above family is
obvious, at least when $N$ becomes large. To achieve an asymptotic
equality between left and right sides in (\ref{nomerr}) we
choose $\t=\t_N \nearrow 1$ as $N\to\iy.$ For instance, if we
choose $\t_N=1-\f cN,\ c>0$, we get
$$\|k_N(\cdot;\t)\|_2^2=\f N\pi+O(1)
$$
and, at the same time,
$$\|k_N(\cdot;\t)\|_1=\|k(\cdot;\t)\|_1=O(\log N).$$
\vn\textbf{Definition.} We call a sequence of kernels $k_N(x)\in\Mc_N$ an
\textit{asymptotic} $(\Ec_N,\Lb_1)$-{\it multiplier,} if \nl

\vspace{0.3cm}

 \n a) $\h k_N(t)\equiv 1$ for $|t|\le N$,

 \nl \n b) $\f\pi
N\|k_N\|_2^2\to 1$, as $N\to\iy$,

\nl \n c) $\|k_N\|_1=O(\log N)$, as $N\to\iy$.

\vn \textbf{Remark 3.} \textit{Our
definition catches some of the characteristic features of
multipliers, as they are defined in Approximation Theory, but is
somewhat different by taking into account the stochastic nature of
nonparametric estimation.}
\textit{Neither of the classical kernels are
multipliers, in the above sense. The {\it sinc kernel} satisfies a)
and b), with an exact equality $\|k\|_2^2=N/\pi,$ but is not in
$\Lb_1$. The Fej\'er kernel does not satisfy a)  whereas the
Vall\'ee-Poussin kernel does not satisfy b).} 

\subsection{Definition of the kernel-type estimator}
Following the discussion of the previous subsection, we attempt to estimate $f$ by means of the
kernel-type estimator
 \begin{gather}\label{tif}
 f_n(x)=\frac{1}{n}\sum_{i=1}^n
 k_{h_n}(x-X_i),\quad x\in \Rb,
 \end{gather}
where  $
k_{h_n}(x)={h_n}^{-1}k_n\left(h_n^{-1}{x}\right)$
is the scaled Fej\'{e}r-type kernel based on (see Section 3.1)
\begin{gather}\label{Fk}
k_n(x)=\frac{\cos(\t_n x)-\cos(x)}{\pi(1-\t_n)x^2},
\end{gather}
with $$\quad h_n=\frac{\t_n}{N},\quad \t_n=1-\frac{1}{N},\quad N=\frac{\log
n}{2\g}.$$

\medskip
\n A remarkable feature of our estimator $f_n(x)$ is that
the sequence $k_{h_n}(x)$ is an asymptotic $({\cal E}_N,\Lb_1)$-multiplier (with $N=\frac{\log n}{2\g}$).
To see this, observe that (see Section 3.1)

\begin{gather}\label{k1}
k_n(x)=k_n(-x),\quad \int_{\Rb}k_n(x) dx=1,\quad
\max\limits_{x\in\Rb}k_n(x)=k_n(0)=\dfrac{1}{\pi};
\end{gather}
\begin{gather}\label{k2}
\hat{k}_n(t)=1,\;\mbox{ for}\;|t|\leq \t_n;
\end{gather}
\begin{gather}\label{k3}
\|k_n\|_2^2=\dfrac{1+\t_n}{2\pi}=\dfrac{1}{\pi}+O\left(N^{-1}\right),\quad
n\to\iy;
\end{gather}
\begin{gather*}
\dfrac{4}{\pi^2}\log\dfrac{1+\t_n}{1-\t_n}+\dfrac13\leq
\|k_n\|_1\leq \dfrac{4}{\pi^2}\log\dfrac{1+\t_n}{1-\t_n}+2,
\end{gather*}
and hence
\begin{gather}\label{k4}
\|k_n\|_1=O(\log N),\quad n\to\iy.
\end{gather}
Therefore
$$\hat{k}_{h_n}(t)\equiv 1\quad\mbox{ for }\; |t|\le N,$$ and as $n\to \iy$
\begin{gather}\label{plpl}
\f\pi N\|k_{h_n}\|_2^2\to 1,\quad
\|k_{h_n}\|_1=O(\log N).
\end{gather}
In subsequent discussion,  the fact that $k_{h_n}(x)$ is an asymptotic $({\cal E}_N,\Lb_1)$-multiplier will allow us to apply inequality (\ref{eqn2})
for establishing some useful properties of this function. 

 Note also that the parameters $h_n,$ $ \t_n$ and $N$ are chosen to satisfy

\begin{gather*}
h_n\to 0\quad\mbox{and} \quad n h_n\to \iy\quad \mbox{as}\; n\to \iy.
\end{gather*}

\subsection{Some properties of the scaled kernel}
In this section, we establish further asymptotic properties of the scaled kernel
$$k_{h_n}(x)={h_n}^{-1}k_n\left(h_n^{-1}x\right).$$

We have
\begin{gather}\label{fmax0}
\max\limits_{x} k_{h_n}(x)=k_{h_n}(0)=\frac{1}{\pi h_n},
\end{gather}

\n and, uniformly in $f\in\mathbb{F}(\g,p,M)$, as $n\to \iy$

\begin{gather}\label{ef0}
\E_f k_{h_n}(x-X)=f(x)+O(n^{-1/2}\log N ),
\end{gather}
\begin{gather}\label{varf0}
\E_f k^2_{h_n}(x-X)=\frac{f(x)}{\pi h_n}+O(\log N).
\end{gather}

 Equality (\ref{fmax0}) is obvious. Let us prove (\ref{ef0})  and (\ref{varf0}).
We get from (\ref{eqn2}), (\ref{k4}), and (\ref{fmax0}) that, uniformly over
$\mathbb{F}(\g, p, M)$,

\begin{eqnarray*}
\sup_{x\in\Rb}|\E_f k_{h_n}(x-X)-f(x)|&=&\|k_{h_n}\ast
f-f\|_{\iy}\leq \frac{8M(1+\|k_{h_n}\|_1)}{\pi}e^{-\g
N}\\ &=&O(n^{-1/2}\log N),
\end{eqnarray*}
which yields (\ref{ef0}).

\medskip
 The proof of (\ref{varf0}) is a more delicate problem. For given $f\in{\mathbb{F}}(\g,p,M)$ and
$x\in\Rb$, consider the function

$$k_{h_n}(x-t)f(t),\quad t\in\Rb,$$

\n and observe that
\begin{gather*}
\E_f k_{h_n}^2(x-X)= \int_{\Rb}k^2_{h_n}(x-t)
f(t)dt=\left(k_{h_n}\ast\left[ k_{h_n}(x-\cdot)f(\cdot)\right]\right)(x).
\end{gather*}
In order to prove (\ref{varf0}) we shall apply inequality (\ref{eqn2}).
We can do this because $k_{h_n}(x)$ is an asymptotic $({\cal E}_N,\Lb_1)$-multiplier  (see Section 3.2 for details),
and for any $f\in{\mathbb{F}}(\g,p,M)$ and all
$x\in\Rb$, the function $k_{h_n}(x-t)f(t)$ belongs (as a function of $t$) to the class ${\mathbb{F}}(\g,p,M_{n}),$ with a
constant $M_n=O(n^{1/2})$, as $n\to\iy.$

The latter fact is easy to verify. Indeed, for any $x\in\Rb$,
$\g>0,$ and all $1\leq p\leq \iy$, uniformly over
$\mathbb{F}(\g, p, M)$,

\begin{gather*}
\|{\rm Re}\,
k_{h_n}(x-(\cdot\pm i\g))f(\cdot\pm i\g)\|_p\leq M
\sup\limits_{t\in\Rb}|k_{h_n}(x-(t\pm i\g))|\\
=\frac{M N h_n}{\pi}\sup\limits_{t\in\Rb}\frac{\left|
\cos\left(Nx-N(t\pm i\g)\right)-\cos\left(h_n^{-1}x-h_n^{-1}(t\pm
i\g)\right)\right| }{|x-(t\pm i\g)|^2}\\ \leq M_{\g} e^{\g
N}=M_{\g}n^{1/2}.
\end{gather*}

\n Therefore, noting that when $t=x$ one has

$$k_{h_n}(0)f(x)=\frac{f(x)}{\pi h_n},$$

\medskip
\n we infer from (\ref{eqn2}) and (\ref{plpl}) that, uniformly over
$\mathbb{F}(\g, p, M)$,
\begin{gather*}
\sup\limits_{x\in\Rb}\left|\E_f k_{h_n}^2(x-X)-(\pi
h_n)^{-1}{f(x)}\right|\\ =\sup\limits_{x\in\Rb} \left|\left(k_{h_n}\ast
\left[k_{h_n}(x-\cdot)f(\cdot)\right]\right)(x)-k_{h_n}(x-x)f(x)\right| \\ \leq
 \frac{M_{\g}n^{1/2}(1+\|k_{h_n}\|_1)}{\pi}\, e^{-\g N}=O(\log
 N),
 \end{gather*}
which yields (\ref{varf0}).

\section{Main result}
The main result of the paper deals with an upper bound on the local
maximal $L_p$-risk of our estimator $f_n$ of
density $f\in \mathbb{F}(\g,p,M)$, and provides a partial solution to the conjecture  of
 Guerre and Tsybakov (see relation (\ref{GT})).
Before stating the result, for $1\leq p<\iy$ put
\begin{gather}\label{betam}
\quad\b_p=\pi^{-1/2}\|f\|_{p/2}^{1/2}M_p,\quad M_p=\left(\E|{\cal N}(0,1)|^p\right)^{1/p}=\sqrt{2}\left(\frac{1}{\sqrt{\pi}}\Gamma\left(\frac{p+1}{2}\right)\right)^{1/p},
\end{gather}
and define the \textit{rate function} by
\begin{gather}\label{ratef}
\psi_p(n)={(n h_n)}^{-1/2}\beta_p.
\end{gather}
Denote by ${\cal L}={\cal L}(A,B)$  a class of loss
functions that consists of non-decreasing functions $l:[0,\iy)\to\Rb$ such that $l(0)=0$, $l(x)$ is continuous at $x=1$, and for
some positive constants $A$ and $B$,
\begin{gather}\label{omega}
l(x)\leq A e^{B |x|},\quad x\in[0,\iy). \end{gather}
We have the following theorem.

\bigskip
\n \textbf{Theorem 1.} {\it Let $X_1, X_2,\ldots$ be a sequence of i.i.d. random variables with density function $f_0\in\mathbb{F}(\g,p,M)$,
and  let $f_n$ be an estimator of $f_0$ given by {\rm
(\ref{tif})}. Then
for any $l\in{\cal L}$ and any} $1\leq p< \iy,$

 \begin{eqnarray*}
\lim\limits_{{\cal V}\searrow f_0}
\limsup\limits_{n\to\iy}\sup_{f\in{\cal V}}\Eb_f l
(\psi_p^{-1}(n)\|{f}_n-f\|_p)\leq l(1),
\end{eqnarray*}
\medskip

 \n {\it where ${\cal V}$ is an arbitrary vicinity of $f_0$ in the topology $\mathbb{T}(\g,p,M)$.}

\section{Proof of Theorem 1}
Let $f\in\mathbb{F}(\g,p,M)$ and let $f_n(x)$ be an estimator of $f(x)$ given by (\ref{tif}). Consider the variance-bias decomposition
\medskip
\begin{eqnarray}\label{decomp}
{f}_n(x)-f(x)&=&({f}_n(x)-\Eb_f {f}_n(x))+(\Eb_f {f}_n(x)-f(x))\nonumber\\
& &\nonumber\\
&=&{(n h_n)}^{-1/2}\xi_n(x)+b_n(x),\quad x\in\Rb,
\end{eqnarray}
\medskip
where $\xi_n(x)=\xi_n(x,X_1,\ldots,X_n)$ is an infinitely
differentiable zero-mean stochastic term and $b_n(x)$ is a bias term
defined as follows:
\medskip
\begin{eqnarray}
\xi_n(x)&=&\left(\frac{h_n}{n}\right)^{1/2}\sum\limits_{i=1}^{n}\left(k_{h_n}(x-X_i)-
\Eb_f k_{h_n}(x-X_1)\right),\nonumber \\
 \nonumber \\b_n(x)&=&\Eb_f k_{h_n}(x-X_1)-f(x).\nonumber
\end{eqnarray}
The proof of Theorem 1 is based on the following four lemmas proved in Section 7.

\bigskip

\n \textbf{Lemma 1.} {\it For any $1\leq p<\iy,$ locally uniformly in
$f\in{\mathbb{F}}(\g,p,M)$,

$$\E_f\|\xi_n\|_p^p\to
\b_p^p,\quad n\to \iy,$$

\n where $\beta_p$ is given by {\rm (\ref{betam})}.}
\bigskip

\n \textbf{ Lemma 2.} {\it For any $1\leq p<\iy,$ locally uniformly in
$f\in{\mathbb{F}}(\g,p,M)$,}

$$\Var_f\|\xi_n\|_p^p\to 0,\quad n\to \iy.$$

\bigskip
\n The next lemma deals with the bias term $b_n$.
\bigskip

 \n \textbf{Lemma 3.} {\it For
any $1\leq p<\iy$, locally uniformly in
$f\in{\mathbb{F}}(\g,p,M)$,}

\begin{gather*}\label{bias}
\|b_n\|_p=o\left(\E_f\|{f}_n-\E_f {f}_n\|_p\right).
\end{gather*}

\bigskip
\n The fourth lemma ensures the uniform integrability of the
sequence $\{l(\psi_p^{-1}(n)\|{f}_n-f\|_p)\}$ uniformly in a vicinity of
$f_0\in{\mathbb{F}}(\g,p,M)$.

\bigskip
 \n \textbf{Lemma 4.} {\it For any  $l\in{\cal L}$ and $1\leq p< \iy,$
there exists $\tau>0$ such that, locally uniformly in
$f\in{\mathbb{F}}(\g,p,M)$,}

\begin{gather*}
\limsup_{n\to\iy} \E_f
l ^{1+\tau}(\psi_p^{-1}(n)\|{f}_n-f\|_p)<\iy.
\end{gather*}

\bigskip
\textbf{Proof of Theorem 1.} Lemmas 1 and 2 in conjunction with Chebyshev's inequality imply, locally uniformly in $f\in\mathbb{F}(\g,p,M)$,

\begin{gather*}
\psi^{-1}_p(n)\|{f}_n-\E_f{f}_n\|_p\;{\bd  P\over \to}\;1,
\end{gather*}
\medskip

\n
Since $l(x)$ is continuous at $x=1$, the combination of this relation with
Lemmas 3 and 4 yields via (\ref{decomp}) the conclusion of Theorem 1. \done

\section{Proofs of Lemmas}
We shall first formulate several preliminary results that will be used in the proofs of Lemmas 1--4.
\subsection{Propositions}
In the next section, under the conditions of Theorem 1, we shall prove three technical propositions
connected to the stochastic error $\xi_n(x)$.
\bigskip

\n {\textbf{Proposition 1.}} \textit{Let $\xi_n(x)$ be the stochastic term in decomposition {\rm (\ref{decomp})}, where $X_1,X_2\ldots,$
is a sequence of i.i.d. random variables  with density function $f\in{\mathbb{F}}(\g,p,M)$.} \textit{For any }$p>2$ \textit{we have for
some constant} $L_{p}>0$, \textit{and all} $r\geq1$
\textit{and} $n\geq1$,%

\begin{equation}
\E_f\left\Vert \xi_{n}\right\Vert _{p}^{r}\leq\left(\frac{r
L_{p}}{1+\log r}\right)^r\left[\frac{1}{(n h_n)^{r/2}}\left( \left(  \int_{\Rb}%
f^{p/2}\left(  y\right)  dy\right) ^{r/p}+\frac{1}{\left(
nh_{n}\right)^{r/2-r/p}}\right)  +\frac{h_{n} ^{r/p}}{\left(
nh_{n}\right)^{r}}\right]  .
\label{DD66}%
\end{equation}

\medskip
\n \textit{For }$p=2$ \textit{we have for some constant} $L_{2}>0$,
\textit{and
all} $r\geq1$ \textit{and} $n\geq1$,%

\begin{equation}
\E_f\left\Vert \xi_{n}\right\Vert _{2}^{r}\leq \left(\frac{r
L_{2}}{1+\log r}\right)^r\left[\frac{1}{(n h_n)^{r/2}}+\frac{h_{n}^{r/2}}{(n
h_n)^{r}}\right]  . \label{DD77}%
\end{equation}

\medskip
\n \textit{Moreover, for a given }$1\leq p<2$ {\it and some
${(2-p)}/{p}<\l\leq 2$},
{\it we have for some constant} $L_{p}>0$ \textit{and} \textit{all
}$r\geq1$
\textit{and} $n\geq1$,%

\begin{gather}
\E_f\|\xi_n\|_p^r  \leq \left(\frac{r L_{p}}{1+\log
r}\right)^r\left[\frac{\left(C^{1/p}(\l,p)\left( 1+2^{\lambda-1}\left(
\E|X|^{\lambda}+h_n^{\lambda }\E|Y_n|^{\lambda} \right)\right)
^{1/2}\right) ^{r}}{(n h_n)^{r/2}}+\frac{h_{n}^{r/p}}{\left(
nh_{n}\right)  ^{r}}\right], \label{DD7}%
\end{gather}
\textit{where} $X$ {\it has density $f(x)$ and $Y_n$ has density}
$k_n^2(x)/\|k_n\|_2^2$.

\bigskip
An immediate consequence of Proposition 1 is the following corollary.

\bigskip
\n \textbf{Corollary 1.} \textit{ For any  $1\leq p<\iy$ and $t>0$,}
\textit{whenever the conditions of Theorem {\rm{1}} hold}, \textit{locally uniformly in
$f\in{\mathbb{F}}(\g,p,M)$,}
\begin{gather*}
\limsup_{n\to\iy} \E_f \exp\left(t\sqrt{n
h_n}\|f_n-\E_f f_n\|_p\right)<\iy.
\end{gather*}

\n In terms of the process $\xi_n(x)$ we have for every $t>0$, locally uniformly in
$f\in{\mathbb{F}}(\g,p,M),$
\begin{gather*}
\limsup_{n\to\iy} \E_f
\exp\left(t\|\xi_n\|_p\right)<\iy,
\end{gather*}

\bigskip
\n \textbf{Proof of Corollary 1.} Observe that for $2\leq p<\iy$ the norm $\|f\|_{p/2}$ of $f\in{\cal V}$ can be made as close to $\|f_0\|_{p/2}$ as desired,
since a vicinity ${\cal V}$ of $f_0$ can be made arbitrarily small. Similarly, for $1\leq p<2$
 the expectation $\E_f|X|^{\l}$ can be made as close to $\E_{f_0}|X|^{\l}$ as desired.
Also, for all $1\leq p<2$ and all $(2-p)/p<\l\leq 2$, as $n\to \iy$

\begin{gather}
h_n^{\l}\,
\E|Y_{n}|^{\l}=\frac{8 h_n^{\l}}{\pi^2 \|k_n\|_2^2 }\int_{0}^{\iy}\frac{
\sin^2\left(y(1+\t_n)/2\right)\sin^2\left(y\left(1-\t_n\right)/2\right)}{(1-\t_n)^2
y^{4-\l}} \,dy \nonumber\\
=\frac{8 h_n^{\l}}{\pi^2 \|k_n\|_2^2 (1-\t_n)^2 }\left(\int_0^1 +\int_1^{h_n^{-1}}+\int_{h_n^{-1}}^{\iy}\right)
\frac{\sin^2\left(y(1+\t_n)/2\right)\sin^2\left(y\left(1-\t_n\right)/2\right)}{y^{4-\l}} \,dy\nonumber \\
=O\left( h_n^{\l-2}\right)\int_0^1 \frac{(y(1+\t_n)/2)^2 (y(1-\t_n)/2)^2 }{y^{4-\l}}\, dy +
O\left( h_n^{\l-2}\right)\int_1^{h_n^{-1}} \frac{(y(1-\t_n)/2)^2}{y^{4-\l}}\, dy\nonumber \\ + O\left( h_n^{\l-2}\right)\int_{h_n^{-1}}^{\iy}\frac{dy}{y^{4-\l}}
=O(h_n^{\l})+O(h_n)+O(h_n)=o\left(1\right). \label{eqyn}
\end{gather}
Therefore, 
whenever the
conditions of Theorem 1 hold,
for any $1\leq p<\iy$ and any sufficiently small vicinity ${\cal V}$ of $f_0$, we have, using Proposition 1, that for some constant
$A>0$, and for all $r\geq 1$ and all sufficiently large $n$, 
\begin{gather*}
(nh_n)^{r/2}\sup_{f\in{\cal V}}\E_f\|f_n-\E_f f_n\|_p^r\leq
\left(\frac{r A}{1+\log r}\right)^r\left(1+\frac{1}{(n
h_n)^{r/2}}\right).
\end{gather*}
From this using Stirling's approximation, which says for $r\geq 1$
that $r!>(r/e)^r$, and the assumption $n h_n\to \iy$, for
sufficiently large $n$,
\begin{gather*}
\sup_{f\in{\cal V}}\E_f\exp\left(t(n h_n)^{1/2}\|f_n-\E_f
f_n\|_p\right)\\
\leq\sum_{r=0}^{\iy}\frac{(n h_n)^{r/2}t^r}{ r!}\sup_{f\in{\cal V}}\E_f\|f_n-\E_f f_n\|_p^r\leq 1+
\sum_{r=1}^{\iy}\frac{(tAr)^r }{(1+\log r)^r r!}\left(1+\frac{1}{(n
h_n)^{r/2}}\right)
 \\ \leq 1+
\sum_{r=1}^{\iy}\frac{(t Ae)^r}{(1+\log r)^r} <\iy.
\end{gather*}
\done
\bigskip

\n {\textbf{Proposition 2.}} \textit{Under the conditions of
Theorem 1}, {\it  for all $p\geq 1$ and $r\geq 1$, for any $\e>0$  and any sufficiently small vicinity ${\cal V}$ of $f_0$, there exist
$n^{*}>0$ and $B>0$
such that  for all $n\geq n^{*}$ and all $f\in{\cal V}$}

$$\E_f \left(\int_{|x|>B}|\xi_n(x)|^p\, dx\right)^r< \e.$$

\bigskip
In order to state the next proposition,
we need a concept of the \textit{uniform weak convergence} (see, for example, Ibragimov and Hasminskii (1981), p. 365.)

\medskip
\n \textbf{Definition.}  Let $\eta_{\t}$ and $\eta_{n\t}$, $n=1,2,\ldots,$ be random $k$-vectors with respective distributions ${ \Pb}_{\t}$
and ${\Pb}_{n\t}$, $n=1,2,\ldots$, that depend on a parameter $\t\in\Theta$. We say that
${\eta}_{n\t}$ \textit{converges weakly to} ${\eta}_{\t}$  \textit{uniformly in $\t \in\Theta$}, if for any continuous bounded function $g:\Rb^k\to \Rb$,
uniformly over $\Theta$, as $n\to \iy$
$$\int_{\Rb^k}g(x)\, d{\Pb}_{n\t}(x) \to \int_{\Rb^k}g(x) \,d{\Pb}_{\t}(x) \quad \mbox{or}\quad
\E_\t g(\eta_{n\t})\to \E_\t g(\eta_\t).$$

In what follows, we think of vectors in $\Rb^k$ as row vectors.
 The next statement is a  uniform version of the Central Limit Theorem (CLT) whose proof is similar to that of Theorem 15 in Appendix I of Ibragimov and Hasminskii (1981).
\medskip

\n \textbf{Fact 1.} \textit{Let $\xi_{1,n},\ldots,\xi_{n,n}$, $n=1,2,\ldots,$ be a sequence  of series of random  $k$-vectors
that are independent within each series, and such that the distribution ${\Pb}^{\t}_{i,n}$ of $\xi_{i,n}$, $i=1,\ldots,n$, depends on a parameter $\t\in\Theta$.
Let $\E_{\t} \xi_{i,n}={0}\in\Rb^k$, $\zeta_n=\sum_{i=1}^n \xi_{i,n},$ and
 $\E_{\t}|\xi_{i,n}|^2<\iy$, where $|x|^2=(x,x)$ is the square of the norm of $x$ . Put
$$\sigma_n^2(\theta)= \sum_{i=1}^n\E_{\t} \xi_{i,n}^{\top}\xi_{i,n}.$$
If $\sigma_n^2(\t)\to \s^2(\t)$ uniformly in $\t\in\Theta$
and the Lindeberg condition holds uniformly over $\Theta$, that is,
for any $\tau>0$
\begin{gather}\label{Lind}\sup_{\t\in\Theta} \sum_{i=1}^n \E_{\t}\left(|\xi_{i,n}|^2 \mathbb{I}\left(|\xi_{i,n}|>\tau\right)\right)\to 0,\quad n\to\iy,
\end{gather}
then uniformly in $\t \in\Theta$}
$${\zeta_n} \stackrel{{d}}{\longrightarrow}{N}\left({0},{ \sigma^2(\t)}\right) ,\quad n\to\iy.$$

\medskip
We will apply the preceding fact to the sequence $(\xi_n(x), \xi_n(y))$, $n=1,2,\ldots,$
where $\xi_n(x)$
is the stochastic error of $f_n(x)$ in the decomposition (\ref{decomp}) and $x\neq y$.
To this end, for $n=1,2,\ldots$ and $i=1\ldots,n,$ put
\begin{gather}\label{etab}
{\etab}_{i,n}(x,y)=({\eta}_{i,n}(x),{ \eta}_{i,n}(y)),
\end{gather}
where
$$\eta_{i,n}(x)=\sqrt{\frac{h_n}{n}}\left[{k_{h_n}(x-X_i)-\E_f{k}_{h_n}(x-X_1)}\right],\quad x\in \Rb,$$
and observe that
$$(\xi_n(x), \xi_n(y))=\sum_{i=1}^n {\etab}_{i,n}(x,y),\quad n=1,2,\ldots.$$
Now fix $\delta\in(0,1)$ and  consider the region

\begin{gather}\label{regD}
D_n=D_n(\delta)=\{(x,y)\in\Rb^2: |x-y|\geq  N^{-(1-\delta)/2}\},\quad N=\frac{\log n}{2\g}.
\end{gather}
When $(x,y)\in D_n,$ uniformly in $f\in{\mathbb{F}}(\g,p,M)$, the components of $(\xi_n(x),\xi_n(y)$ are weakly correlated and for all large enough $n$ the uniform CLT applies.

 \bigskip

\n \textbf{Proposition 3.} \textit{Under the conditions of Theorem 1, for all $p\geq 1$ and
any $\e>0$ there exists  $n^*>0$ such that
 for all $n\geq n^*$, any continuous bounded function $g:\Rb^2\to \Rb$, and all $f\in\mathbb{F}(\g,p,M)$,
 \begin{gather}\label{clt2}
\sup_{(x,y)\in D_n}\left|\E_f g((\xi_n(x),\xi_n(y)))- \E_f g({\zetb}(x,y))\right|<\e,
\end{gather}
where ${\zetb}(x,y)$ is a mean zero normal random vector with covariance matrix $\Cov_f(\zetb(x,y))={\rm Diag}\left(\pi^{-1}{f(x)},\pi^{-1}{f(y)}\right).$}

\subsection{Proofs of Lemmas 1--4}
The proofs of Lemmas 1--4 are largely based on Propositions 1--3.

\bigskip
\n \textbf{Proof of Lemma 1.} First, we show that for any $p\geq 1$,
uniformly over a small vicinity ${\cal V}$ of $f_0$, as $n\to \iy$
\begin{gather}\label{lem33}
\E_f|\xi_n(x)|^p\to\E\left|\sqrt{\frac{f(x)}{\pi}}\xi\right|^p=\frac{f^{p/2}(x)}{\pi^{p/2}}M_p^p,\quad
 x\in\Rb,
\end{gather}
where $\xi$ is a standard normal random variable and $M_p=\left(\E|{\cal N}(0,1)|^p\right)^{1/p}$. Indeed, by the CLT
and the continuity theorem, as $n\to \iy$
\begin{gather*}
|\xi_n(x)|^p\;\;{\bd d\over
\to}\;\;\frac{f^{p/2}(x)}{\pi^{p/2}}|\xi|^p,\quad x\in\Rb.
\end{gather*}
Therefore, in order to prove (\ref{lem33}) it is sufficient to show that  there exists $\tau>0$ such that
\begin{gather}\label{integr}
\lim_{{\cal V}\searrow f_0}\limsup_{n\to\iy}\sup_{f\in{\cal V}}\E_f\left(|\xi_n(x)|^p\right)^{1+\tau}<c<\iy.
\end{gather}
Note that for any $p>0$,
\begin{gather*}
|x|^p\leq ce^{|x|}<c(e^{x}+e^{-x}),\quad x\in\Rb.
\end{gather*}
Thus (\ref{lem33}) holds provided for a sufficiently small vicinity ${\cal V}$ of $f_0$
\begin{gather}\label{plmin0}
\limsup_{n\to\iy}\sup_{f\in{\cal V}}\E_f e^{\pm\xi_n(x)}<c<\iy.
\end{gather}
The proofs of both cases in (\ref{plmin0}) are similar. Let us
consider the case of the plus sign. By independence
\begin{gather*}
\E_f \exp{\{\xi_n(x)\}}
=\left(\E_f\exp\left\{\sqrt{\frac{h_n}{n}}\left(k_{h_n}(x-X_1)-\E_f
k_{h_n}(x-X_1)\right)\right\}\right)^n.
\end{gather*}
Using Tailor's expansion, formulas (\ref{fmax0})--(\ref{varf0}), and the fact that $f\in\mathbb{F}(\g,p,M)$ is a bounded function,
for some $C>0$ and all sufficiently large $n$, uniformly in $f\in\mathbb{F}(\g,p,M)$,
\begin{gather*}
\E_f\exp\left\{\sqrt{\frac{h_n}{n}}(k_{h_n}(x-X_1)-\E_f
k_{h_n}(x-X_1))\right\}\\
=\E_f\left(\sum_{j=0}^{\iy}\frac{1}{j!}\left(\sqrt{\frac{h_n}{n}}\left(k_{h_n}(x-X_1)-\E_f
k_{h_n}(x-X_1)\right)\right)^j\right)\\= 1+\frac{h_n}{2n}{\Var
k_{h_n}(x-X_1)}+\E_f\left(\sum_{j=3}^{\iy}\frac{1}{j!}\left(\sqrt{\frac{h_n}{n}}(k_{h_n}(x-X_1)-\E_f
k_{h_n}(x-X_1))\right)^j\right)\\ \leq
1+\frac{C}{2n}+\sum_{j=3}^{\iy}\frac{(C\log n)^j}{n^{j/2}(\log
n)^{j/2}j!}=1+\frac{C}{2n}+\frac1n\sum_{j=3}^{\iy}\frac{(C\log
n)^{j/2}}{n^{(j-2)/2}j!}\leq 1+\frac{C}{n}.
\end{gather*}
Therefore for any $p>0$, uniformly in $f\in\mathbb{F}(\g,p,M),$
\begin{gather}\label{otsen0}
\limsup_{n\to\iy}
\E_f|\xi_n(x)|^p< (2c) \limsup_{n\to\iy} \left(1+\frac{C}{n}\right)^n=\mbox{const}<\iy,
\end{gather}
and relation (\ref{plmin0}), and hence relation (\ref{integr}), follows.

Next, by conditions {\bf (A1)} and {\bf(A2)}, Remark 2, and Proposition 2, for any $\e>0$ and any sufficiently small vicinity ${\cal V}$ of $f_0$,
there exist numbers $n^*>0$ and $B>0$ such that for all $n\geq n^*$ and all $f\in{\cal V}$ the inequalities
$$\sup_{f\in {\cal V}}\int_{|x|>B} f^{p/2}(x) dx<\e\quad\mbox{and}\quad \sup_{f\in {\cal V}}\int_{|x|>B}\E_f|\xi_n(x)|^p dx<\e$$
hold  simultaneously.
Thus, we only need to show that, locally uniformly in $f\in\mathbb{F}(\g,p,M)$,
\begin{gather*}
\int_{-B}^B\E_f|\xi_n(x)|^p
dx\to\int_{-B}^B\E_f\left|\sqrt{\frac{f(x)}{\pi }}\xi\right|^p
dx=
\beta_p^p,\quad n\to \iy.
 \end{gather*}
Using the Dominated Convergence Theorem, we get from (\ref{lem33}) and (\ref{otsen0}) that  the above
 limiting relation  holds true.
The proof of Lemma
1 is completed. \done
\bigskip

\n \textbf{Proof of Lemma 2.} In the proof, the key role
belongs to Proposition 2. From Proposition 2, for any $\e>0$ and any sufficiently small
vicinity ${\cal V}$ of $f_0$, there exist $n^{*}>0$ and $B>0$ such
that for all $n\geq n^{*}$,

$$\sup\limits_{f\in{\cal V}}\left|\Var_f\left( \int_{\Rb}\xi_n(x)|^pdx\right)-\Var_f\left(
\int_{|x|\leq B}|\xi_n(x)|^pdx\right)\right|< \e.$$
Due to the identity
\begin{gather}\label{f}
\Var_f \left(\int_{|x|\leq B}|\xi_n(x)|^p
dx\right)=\iint_{[-B,B]^2}\Cov_f\left(|\xi_n(x)|^p,|\xi_n(y)|^p\right)dx dy,
\end{gather}
the problem is reduced to showing that
$$\lim_{{\cal V}\searrow f_0}\lim_{n\to\iy}\sup_{f\in{\cal V}}\iint_{[-B,B]^2}\Cov_f\left(|\xi_n(x)|^p,|\xi_n(y)|^p\right)dx dy=0.$$

As in Proposition 3, consider the region

$$D_n=\{(x,y)\in\Rb^2: |x-y|\geq  N^{-(1-\delta)/2}\},\quad N=\frac{\log n}{2\g},$$

\n with $\delta\in(0,1)$ being fixed, and let $D^c_n$ be its complement (in $\Rb^2$).
\bigskip
\n Since
$$\left|[-B,B]^2\cap D^c_n\right|=O(N^{-(1-\delta)/2}),$$

\medskip
\n it follows from (\ref{f}) and Proposition 1 that

$$\lim_{{\cal V}\searrow f_0}\lim_{n\to\iy}\sup_{f\in{\cal V}}\iint_{[-B,B]^2\cap D^c_n}\Cov_f\left(|\xi_n(x)|^p,|\xi_n(y)|^p\right) dx dy=0,$$
and it remains to show that

\vspace{-0.5cm}
\begin{gather}\label{covar}
\lim_{{\cal V}\searrow f_0}\lim_{n\to\iy}\sup_{f\in{\cal V}}
\iint_{[-B,B]^2\cap D_n}\Cov_f\left(|\xi_n(x)|^p,|\xi_n(y)|^p\right) dx dy= 0.
\end{gather}

In order to prove (\ref{covar}) we shall use Proposition 3 and the fact that for $1\leq p<\iy$, uniformly in $f\in\mathbb{F}(\g,p,M)$,
the sequences
$\{|\xi_n(x)|^p \}$ and $\{|\xi_n(x)\xi_n(y)|^p\}$ are uniformly integrable.
The uniform integrability of
$\{|\xi_n(x)|^p \}$ has been already verified (see (\ref{otsen0})). The sequence
$\{|\xi_n(x)\xi_n(y)|^p\}$ is treated similarly. For any $p>0$ and all $x,y\in
\Rb$
\vspace{-0.5cm}

\begin{gather*}
\E_f|\xi_n(x)\xi_n(y)|^p\leq c \E_f
(e^{\xi_n(x)\xi_n(y)}+e^{-\xi_n(x)\xi_n(y)}).
\end{gather*}

\medskip
\n For the plus sign, by independence \vspace{-0.5cm}

 \begin{gather*} \E_f \exp\{\xi_n(x)\xi_n(y)\}\\
=\left(\E_f\exp\left\{ ({n^{-1}}h_n)[k_{h_n}(x-X_1)-\E_f
k_{h_n}(x-X_1)][k_{h_n}(y-X_1)-\E_f k_{h_n}(y-X_1)]\right\}\right)^n.
\end{gather*}

\medskip
\n Then, using (\ref{fmax0})--(\ref{varf0}), 
for some
$C>0$ and all sufficiently large $n$, uniformly in $f\in\mathbb{F}(\g,p,M)$,
\begin{gather*}
\E_f\exp\left\{ ({n^{-1}}h_n)[k_{h_n}(x-X_1)-\E_f
k_{h_n}(x-X_1)][k_{h_n}(y-X_1)-\E_f k_{h_n}(y-X_1)]\right\}\\
=\E_f\left(\sum_{j=0}^{\iy}\frac{1}{j!}\left(({n^{-1}}h_n)[k_{h_n}(x-X_1)-\E_f
k_{h_n}(x-X_1)][k_{h_n}(y-X_1)-\E_f k_{h_n}(y-X_1)]\right)^j\right)\\\leq
1 +\frac{C}{2n}+\sum_{j=2}^{\iy}\frac{(C \log
n)^{2j}}{n^{j}(\log n)^{j} j!}
=1+\frac{C}{2n}+\frac{1}{n}\sum_{j=2}^{\iy}\frac{(C \log
n)^{j}}{n^{j-1} j!}\leq 1+\frac{C}{n}.
\end{gather*}
Hence for any $1\leq p<\iy$, uniformly in $f\in\mathbb{F}(\g,p,M)$,
\begin{gather*}
\limsup_{n\to\iy}\E_f|\xi_n(x)\xi_n(y)|^p<(2c)\limsup_{n\to\iy}\left(1+\frac{C}{n}\right)^n =\mbox{const}<\iy,
\end{gather*}
and the uniform integrability of the sequence $\{|\xi_n(x)\xi_n(y)|^p\}$  follows.

Therefore, by Proposition 3 and the continuity theorem we readily conclude that
\begin{gather*}
\lim_{{\cal V}\searrow f_0}\lim_{n\to\iy}\sup_{f\in{\cal V}}\sup_{(x,y)\in D_n}\Cov_f
(|\xi_n(x)|^p,|\xi_n(y)|^p)= 0,
\end{gather*}

\n which yields (\ref{covar}). The proof is completed. \done
\bigskip

\n \textbf{Proof of Lemma 3.} It was shown in Section 3.2 that the scaled kernel
$$k_{h_n}(x)=h_n^{-1}k_n(h_n^{-1}x),$$
with $k_n(x)$ given by (\ref{Fk}) and $h_n={N}^{-1}\left(1-{N}^{-1}\right),$
 is an asymptotic $({\cal E}_N, \Lb_1)$-multiplier. Therefore it follows
from (\ref{eqn2}) and (\ref{k4}) that for any $1\leq p<\iy$, uniformly in $f\in\mathbb{F}(\g,p,M)$,
\begin{eqnarray*}
\|b_n(\cdot)\|_p&=&\|\E_f k_{h_n}(\cdot-X) -f(\cdot)\|_p=\|k_{h_n}*f(\cdot)-f(\cdot)\|_p\nonumber\\
& &\nonumber\\ &\leq& \frac{8 M
\left(1+\|k_n\|_1\right)}{\pi}\,e^{-\g N}= O(\log N)e^{_-\g
N}=O\left(n^{-1/2}\log N\right).
\end{eqnarray*}

\n At the same time, due to Lemma 1, locally uniformly in $f\in\mathbb{F}(\g,p,M)$,

\vspace{-0.5cm} \begin{gather*}\E_f\|f_n-\E
f_n\|_p=O\left(n^{-1/2}N^{1/2}\right).
\end{gather*}

\medskip
\n Hence, locally uniformly in $f\in\mathbb{F}(\g,p,M)$, \vspace{-0.5cm}

\begin{gather*}\|b_n\|_p=o\left(\E_f\|f_n-\E_f
f_n\|_p\right),\end{gather*}

\medskip
\n and the lemma  is proved. \done
\bigskip

\n \textbf{Proof of Lemma 4.} In view of condition (\ref{omega}),
the lemma is an immediate consequence of Corollary 1. \done.

\section{Proofs of Propositions}

\textbf{Proof of Proposition 1.} Proposition 1 is a
modification of Proposition 2 of Mason (2009). The latter deals
with the kernel estimator
\begin{gather}\label{mest}
\ti{f}_n(x)=\frac{1}{n h_n}\sum_{i=1}^n
K\left(\frac{x-X_i}{h_n}\right),\quad x\in\Rb,
\end{gather}
where $K$ is a fixed kernel function  such that the least decreasing radial
majorant of $K^2/\|K\|^2_2$ defined by $\Psi(x)=\sup\limits_{|y|\geq
|x|}K^2(y)/\|K\|^2_2$, $x\in\Rb$, is integrable. Mason (2009) suggested a
powerful method for deriving good asymptotic bound for
$$\E\int_{A}|\ti{f}_n(x)-\E \ti{f}_n(x)|^p dx,\quad 1\leq p<\iy,$$
with $\ti{f}_n$ given by (\ref{mest}) and $A$ being a measurable
subset of $\Rb$. Moreover, he provided a useful finite sample bound for

$$\E\|\ti{f}_n-\E \ti{f}_n\|_p^r, \quad 1\leq p<\iy,\quad r\geq 1,$$

\n see also Lemma 4 of Ibragimov and Hasminskii (1980).

 Unlike Mason, we allow our kernel $k_n(x)$ to change with
$n$ rather than being fixed. Due to properties
(\ref{k1})--(\ref{k4}), this, however, does not effect the method.
The integrability of
$$\Psi_n(x)=\sup\limits_{|y|\geq
|x|}k_n^2(y)/\|k_n\|^2_2,\quad x\in\Rb,\quad n=1,2,\ldots,$$ where
$k_n$ is the Fej\'{e}r-type kernel given by (\ref{Fk}) is obvious.

When applying Proposition 2 of Mason (2009) to $\xi_n(x)$, we use the fact
that by (\ref{k1}) and (\ref{k3}), for any $p\geq 2$,
\begin{gather}\label{knn}
\|k_n\|_p^p\leq \max_{x}|k_n(x)|^{p-2}\int_\Rb k_n^2(x) dx=k^{p-2}_n(0)\|k_n\|_2^2=
\frac{1}{\pi^{p-2}}\frac{1+\t_n}{2\pi}\leq
\frac{1}{\pi^{p-1}}.
\end{gather}

Note also that for $1\leq p<2$, the method of Mason requires the
existence of a constant $\lambda>{(2-p)}/{p}$ such that
\begin{gather*}
\E|X|^{\lambda}<\iy,\quad \E|Y_n|^{\lambda}<\iy,
\end{gather*}
where the random variables $X$ and $Y_n$ have densities $f(x)$ and
$k_n^2(x)/\|k_n\|_2^2$, respectively. In our case, the existence of
such $\l$ follows from the definition of the class $\mathbb{F}(\g,p,M)$ and the fact that for any $0\leq\l\leq
2$ and $n\geq 1$, cf. Gradshtein and Ryzhik
(2000), 3.828.9, 3.828.10, and (\ref{eqyn}),
\begin{gather}\label{a123}
\E|Y_n|^{\l}=
\frac{8}{\pi^2 \|k_n\|_2^2 }\int_{0}^{\iy}\frac{
\sin^2\left(y(1+\t_n)/2\right)\sin^2\left(y\left(1-\t_n\right)/2\right)}{(1-\t_n)^2
y^{4-\l}} \,dy<\iy.
\end{gather}
When $1\leq p<2$ the set $\{\l:(2-p)/{p}<\l
\leq 2\}$ is non-empty, so that the required constant $\l>(2-p)/p$, for which $\E|Y_n|^\l<\iy$, does exist.
Therefore, for the proof of Proposition 1 we refer to Mason (2009), pp. 75--78.

\bigskip
\n {\textbf{Proof of Proposition 2.}}
In view of conditions {\bf (A1)} and {\bf(A2)}, it suffices to show that
for any
$B>0$, all $p\geq 2$ and $r\geq 1,$
we can find a constant $C_{p,r}$ such that for any function $f\in\mathbb{F}(\g,p,M)$,
\begin{gather}
\limsup_{n\to\iy}\E_f\left(\int_{|x|>B}|\xi_n(x)|^p\,
dx\right)^r\leq C_{p,r} \left(\int_{|x|>B}f^{p/2}(x)\,dx\right)^r,\label{predel1}
\end{gather}
and
for any
$B>0$, all $1\leq p<2$ and $r\geq 1,$
we can find a constant $D_{p,r}$ such that for any function $f\in\mathbb{F}(\g,p,M)$,
\begin{gather}\limsup_{n\to\iy}\E_f\left(\int_{|x|>B}|\xi_n(x)|^p
\,dx\right)^r\leq D_{p,r}\left(\int_{|x|>B}(1+|x|^{\l})f(x)\, dx\right)^r.\label{predel2}
\end{gather}

\n Inequalities (\ref{predel1}) and (\ref{predel2}) will be obtained by modifying the
arguments in the proof of Proposition 1 of Mason (2009). One of the
key elements of that proof is the following
\bigskip

\noindent\textbf{Fact 2. }(\textbf{Theorem 1 of Talagrand
(1989).)} If $\mathbf{B}$ is a separable Banach space with norm
$\left\Vert \cdot \right\Vert $, $Z_{i}$, $i\in \Nb$, are
independent mean zero random vectors taking values in $\mathbf{B}$,
then for a universal constant $D>0$, for all $r\geq 1$ and $n\geq1$,

\begin{equation}
\left(  \E \Vert S_{n}\Vert^{r}  \right)  ^{1/r}\leq
\frac{Dr}{1+\log r}\left(\E\| S_{n}\|+\left(  \E\max_{1\leq i\leq n}\Vert Z_{i}%
\Vert^{r}\right)  ^{1/r}\right), \label{pin33}%
\end{equation}
\medskip
\n where $S_{n}=Z_{1}+\dots+Z_{n}$. Using the $c_{r}$-inequality, we
get from (\ref{pin33}) the bound

\begin{equation}
\E \Vert S_{n}\Vert^{r} \leq \frac{D^r 2^{r-1}r^r}{(1+\log r)^r}
\left(\left( \E\Vert S_{n}\|\right) ^{r}+\E\max_{1\leq i\leq n}\Vert
Z_{i}\Vert^{r}\right).
\label{pin44}%
\end{equation}
\medskip

Adopting Mason's arguments to our needs, we apply the bound (\ref{pin44}) to the random functions%

\begin{gather*}
Z^B_{i}\left(  x\right)  =(n h_n)^{-1}{\left\{k_n\left(  \frac{x\text{ }-X_{i}}{h_{n}%
}\right)  -\E_f k_n\left(  \frac{x\text{
}-X_1}{h_{n}}\right)\right\}\mathbb{I} \{|x|>B\}},\text{ }x\in\Rb,\; i=1,\dots,n,
\end{gather*}
where  $B$ is a positive constant. Because of the properties of $k_n$ (see
Section 3.2), every $Z^B_{i}\left(  x\right)$ is in $\Lb_{p}$,
for any $p\geq1$. The random function defined by
$$S_n^B(x)=(f_n(x)-\E_f f_n(x))\mathbb{I}\{|x|>B\},\quad x\in\Rb,$$ can be now expressed as
follows:
$$S_n^B(x)=Z_1^B(x)+\ldots +Z_n^B(x),\quad x\in\Rb.$$
Using Jensen's inequality,
\begin{eqnarray}\label{BBB0}
\E_f\Vert S^B_{n}\Vert_{p}&=&\E_f\| (f_{n}-\E_f
f_{n})\mathbb{I}\{|\cdot|>B\}\| _{p}=\E_f\left( \int_{|x|>B}\left\vert
f_{n}\left( x\right) -\E_f f_{n}\left( x\right) \right\vert
^{p}dx\right)^{1/p}\nonumber \\ &\leq&\left( \E_f\int_{|x|>B}
|f_{n}(x)-\E
f_{n}(x)|^{p}\,dx\right)^{1/p}=\left(\E_f\|S^B_n\|_p^p\right)^{1/p}.
\end{eqnarray}

Next, according to Mason (2009), pp. 71--72, for any $f\in{\mathbb{F}}(\g,p,M)$ and any $B>0$,
\begin{eqnarray}
\E_f\int_{|x|>B}\!|{f}_n(x)-\E_f {f}_n(x)|^p\, dx&\!\!\!\leq\!\!\!& \frac{1}{n
h_n}\int_{|x|>B}\!\!\!dx\left(\int_{\Rb}\dfrac{1}{h_n}k^2_n\frac{x-y}{h_n}\right)f(y)\,dy,
\quad p=2, \label{m1}\\
\nonumber \\ \nonumber\\ \E_f\int_{|x|>B}\!|{f}_n(x)-\E_f {f}_n(x)|^p\, dx&\!\!\!\leq\!\!\!&
\dfrac{2^p C_p}{(n
h_n)^{p/2}}\int_{|x|>B}\left(\int_{\Rb}\dfrac{1}{h_n}k^2_n\left(\frac{x-y}{h_n}\right)f(y)\,dy\right)^{p/2}dx\nonumber\\
&&+ \dfrac{{2^p}C_p\|k_n\|_p^p}{(n h_n)^{p-1}},\quad
2<p<\iy, \label{m2}\\ \nonumber\\ \nonumber\\ \E_f\int_{|x|>B}\!|{f}_n(x)-\E_f
{f}_n(x)|^p\, dx&\!\!\!\leq\!\!\!& \frac{1}{(n h_n)^{p/2}}\int_{|x|>B}\left(K_{h_n}^2 *
f(x)\right)^{p/2} dx,\quad 1\leq p<2, \label{m3}
 \end{eqnarray}
 where $$K_{h_n}^2(x)=h_n^{-1} k_n^2(h_n^{-1}x),\quad x\in\Rb.$$

\bigskip

\n  \textbf{Case 1a.} $2<p<\iy.$ From (\ref{m2})

\begin{gather*}
\E_f\int_{|x|>B}|f_n(x)-\E_f f_n(x)|^p\, dx \nonumber \\
\leq\frac{2^{p}C_{p}}{(n h_n)^{p/2}}\int_{|x|>B}\left(
\int_{\Rb}\frac{1}{h_n}k_n^{2}\left( \frac{x-y}{h_{n}}\right)
f\left(
y\right)  dy\right)^{p/2}dx+\frac{2^{p}C_{p}\left\Vert k_n\right\Vert _{p}%
^{p}}{\left(  nh_{n}\right)  ^{p-1}}\nonumber
\end{gather*}\vspace{-0.2cm}
\begin{gather}\label{BBB11}
=\frac{2^{p}C_{p}}{(n h_n)^{p/2}}\int_{|x|>B}\left(
\int_{\Rb}\frac{1}{h_n}k_n^{2}\left( \frac{t}{h_{n}}\right) f\left(
x-t\right)  dt\right)^{p/2}dx+\frac{2^{p}C_{p}\left\Vert k_n\right\Vert _{p}%
^{p}}{\left(  nh_{n}\right)  ^{p-1}}.
\end{gather}
Consider  the integral on the right-hand side of (\ref{BBB11}).
Using  the rough version of the $c_r$-inequality that says that for all $r>0$
$$|x+y|^r\leq 2^{r}(|x|^{r}+|y|^{r}),$$
we get
\begin{gather*}
\int_{|x|>B}\left( \int_{\Rb}\frac{1}{h_n}k_n^{2}\left(
\frac{t}{h_{n}}\right) f\left( x-t\right)  dt\right)^{p/2}dx \nonumber\\
=\int_{|x|>B}\left( \int_{|t|\leq B/2}\frac{1}{h_n}k_n^{2}\left(
\frac{t}{h_{n}}\right) f\left( x-t\right)  dt+ \int_{|t|>
B/2}\frac{1}{h_n}k_n^{2}\left( \frac{t}{h_{n}}\right) f\left(
x-t\right) dt \right)^{p/2}dx\nonumber\\
\leq 2^{p/2}\int_{|x|>B}\left( \int_{|t|\leq
B/2}\frac{1}{h_n}k_n^{2}\left( \frac{t}{h_{n}}\right) f\left(
x-t\right)  dt\right)^{p/2}+\nonumber \\
\end{gather*}\vspace{-1.3cm}
\begin{gather}\label{BBB22}
+2^{p/2}\int_{|x|>B}\left( \int_{|t|> B/2}\frac{1}{h_n}k_n^{2}\left(
\frac{t}{h_{n}}\right) f\left( x-t\right)
dt\right)^{p/2}=:{T_1}+T_2.
\end{gather}
The first integral in (\ref{BBB22}) is estimated as follows:

\begin{gather*}
T_1\leq2^{p/2} \int_{|s|>B/2}f^{p/2}(s)\,ds \left(
\int_{\Rb}\frac{1}{h_n}k_n^{2}\left( \frac{t}{h_{n}}\right)
dt\right)^{p/2}\nonumber \\
=2^{p/2}\int_{|s|>B/2}f^{p/2}(s)\,ds\left(\int_{\Rb}k_n^2(x)\,dx\right)^{p/2}=
2^{p/2}\|k_n\|_2^p\int_{|s|>B/2}f^{p/2}(s)\,ds.
\end{gather*}
For the second integral we have

\begin{gather*}
T_2\leq 2^{p/2}\int_{\Rb}f^{p/2}(s)ds \left(
\int_{|t|>B/2}\frac{1}{h_n}k_n^{2}\left( \frac{t}{h_{n}}\right)
dt\right)^{p/2}\\ =2^{p/2}\int_{\Rb}f^{p/2}(s)ds \left(
\int_{|x|>B/2h_n}\!\!\!k_n^{2}\left( x\right)dx\right)^{p/2}.
\end{gather*}
Therefore using (\ref{BBB11}) and (\ref{BBB22})
\begin{gather*}
\E_f\int_{|x|>B}|f_n(x)-\E f_n(x)|^p dx \nonumber \\
\leq\frac{2^{p}C_{p}}{(n
h_n)^{p/2}}\left\{2^{p/2}\|k_n\|_2^p\int_{|x|>B/2}f^{p/2}(x)\,dx+2^{p/2}\int_{\Rb}f^{p/2}(x)\,dx
\left[\int_{|x|>B/2h_n}k_n^{2}\left(
x\right)dx\right]^{p/2}\right\}\nonumber \\+\frac{2^p C_p
\|k_n\|_p^p}{(n h_n)^{p-1}}.
\end{gather*}
From this, due to (\ref{BBB0}) and using the rough version of the $c_r$-inequality, for any
$r\geq 1$,

\begin{gather*}
\left(\E_f\|S_n^B\|_p\right)^r \leq\left\{\E_f\int_{|x|>B}|f_n(x)-\E_f
f_n(x)|^p
dx\right\}^{r/p}
\leq 2^{r/p}\left\{\frac{2^p C_p}{(n
h_n)^{p/2}}\right\}^{r/p}\\ \times\left\{ 2^{p/2} \|k_n\|_2^p
\int_{|x|>B}\!\!\!f^{p/2}(x)\,dx+ 2^{p/2} \int_{\Rb}f^{p/2}(x)\,dx
\left[\int_{|x|>B/2h_n}\!\!\!k_n^{2}\left(
x\right)dx\right]^{p/2}\right\}^{r/p} \\ + 2^{r/p} \left\{\frac{2^p
C_p \|k_n\|_p^p}{(n h_n)^{p-1}}\right\}^{r/p} \leq
\frac{2^{2r/p}2^{3r/2} C_p^{r/p}}{(n h_n)^{r/2}} \|k_n\|_2^{r}
\left(\int_{|x|>B}\!\!\!f^{p/2}(x)\,dx\right)^{r/p}
\end{gather*}
\vspace{-0.6cm}
\begin{gather}\label{snb0}
+\frac{2^{2r/p}2^{3r/2} C_p^{r/p}}{(n
h_n)^{r/2}}\left(\int_{\Rb}f^{p/2}(x)\,dx\right)^{r/p}
\left(\int_{|x|>B/2h_n}k_n^{2}\left( x\right)dx\right)^{r/2} +
\frac{2^{r/p}2^{r} C_p^{r/p}\|k_n\|_p^r}{(n h_n)^{(p-1)r/p}}.
\end{gather}

\medskip
\n Recalling that, as $n\to \iy$,

$$N=O(\log n),\quad\t_n=1-\frac{1}{N}=O(1) ,\quad h_n
=\frac{\t_n}{N}=O(\log ^{-1}n),$$

\medskip \n we get

$$k^2_n(x)=\left|\frac{\cos(\t_n x)-\cos(x)}{\pi (1-\t_n)
x^2}\right|^2\leq \frac{2N^2}{\pi^2 x^4},$$ and hence

\begin{gather}\label{ner11}
\int_{|x|>B/2h_n}k_n^{2}\left( x\right)dx\leq
\frac{2N^2}{\pi^2}\int_{|x|>B/2h_n}\frac{dx}{x^4}=O(N^2
h_n^3)=O(h_n).
\end{gather}

\bigskip

\n Now, by Fact 2 applied to the random process
$S_n^B=Z_1^B+\ldots +Z_n^B$, we have, cf. (\ref{pin44}),

\begin{equation}\label{ner22}
\E_f  \Vert S^B_{n}\Vert_p^{r}  \leq \frac{D^r
2^{r-1}r^r}{(1+\log r)^r} \left(\left( \E_f\Vert S^B_{n}\|_p\right)
^{r}+\E_f\max_{1\leq i\leq n}\Vert Z^B_{i}\Vert_p^{r}\right),
\end{equation}

\medskip
\n where, using  $c_r$-inequality,
\begin{gather}
\E_f\max_{1\leq i\leq n}\|Z_i^B\|_p^r\nonumber \\=
 \E_f\max_{1\leq i\leq n}\left(
\int_{|x|>B}\left\vert \frac {1}{nh_{n}}\left( k_n\left(
\frac{x-X_{i}}{h_{n}}\right)  -\E_f k_n\left( \frac{x-X_1}{h_{n}}\right)
\right)  \right\vert ^{p}dx\right)  ^{r/p}\nonumber \\
 \leq
 2^{r+1} \E\max_{1\leq i\leq n}\left(
\int_{|x|>B}\left\vert \frac {1}{nh_{n}} k_n\left(
\frac{x-X_{i}}{h_{n}}\right)\right\vert ^{p}dx\right) ^{r/p}\leq
\frac{2^{r+1} h_n^{r/p}\|k_n\|_p^r}{(n h_n)^r}.\label{nerav33}
\end{gather}

\medskip
\n We see then by (\ref{snb0})--(\ref{nerav33}) that for all sufficiently large $n$,
\begin{gather*}
\E_f\|S_n^B\|_p^r=\E_f\left(\int_{|x|>B}|f_n(x)-\E_f f_n(x)|^p
dx\right)^{r/p}\\
\leq \frac{D^r 2^{r-1}r^{r}2^{2r/p}2^{3r/2} C_p^{r/p}}{(1+\log r)^r(n
h_n)^{r/2}}
\left\{\|k_n\|_2^{r}\left(\int_{|x|>B/2}f^{p/2}(x)\,dx\right)^{r/p}
+c_\g h_n\int_{\Rb}f^{p/2}(x)\,dx \right\}\\
+\frac{D^r 2^{2r-1}r^{r}2^{r/p} C_p^{r/p}\|k_n\|_p^r}{(1+\log r)^r(n
h_n)^{(p-1)r/p}} +\frac{D^r 2^{2r}r^{r}h_n^{r/p}\|k_n\|_p^r}{(1+\log
r)^r(n h_n)^{r}}.
\end{gather*}
Next, by (\ref{knn}) for $p\geq 2$

$$\|k_n\|_p^p\leq \frac{1}{\pi^{p-1}},$$

\medskip
\n and

$$\frac{(n h_n)^{r/2}}{(n h_n)^{(p-1)r/p}}\to 0,\quad \frac{h_n^{r/p}(n h_n)^{r/2}}{(n h_n)^{r}}\to
0, \quad r\geq 1,\;p> 2.$$

\bigskip \n Therefore
for all $p>2,$ $r\geq
1$, and $B>0$ there exists a number $n^*>0$ such that
for all $n\geq n^*$
\begin{gather*}
\E_f\left(\int_{|x|>B}|f_n(x)-\E_f f_n(x)|^p
dx\right)^{r/p}\\
\leq \frac{a_{p,r}}{(n
h_n)^{r/2}}\left\{\left(\int_{|x|>B/2}f^{p/2}(x)\,dx\right)^{r/p}+
c_\g h_n\int_{\Rb}f^{p/2}(x)\,dx \right\}
\end{gather*}
\vspace{-0.3cm}
\begin{gather*}
+\frac{b_{p,r}}{(n h_n)^{(p-1)r/p}}+ \frac{c_{p,r}\,h_n^{r/p}}{(n
h_n)^r} \leq \frac{C_{p,r}}{(n
h_n)^{r/2}}\left(\int_{|x|>B/2}f^{p/2}(x)\,dx\right)^{r/p},
\end{gather*}
where $a_{p,r}$, $b_{p,r}$, $c_{p,r}$, and $C_{p,r}$ are constants
that may depend on $p$ and $r$.
Noticing that

\begin{gather}\label{plus}
\xi_n(x)=(n h_n)^{1/2}(f_n(x)-\E_f f_n(x))
\end{gather}
yields
(\ref{predel1}).

\medskip
\n \textbf{Case 1b.} $p=2$. In this case, the proof is analogous to
that of Case 1a, but even easier. Using (\ref{m1}) we get
\begin{gather*} \E_f\int_{|x|>B}|f_n(x)-\E_f f_n(x)|^2
dx\leq \frac{1}{n
h_n}\int_{|x|>B}dx\int_{\Rb}\frac{1}{h_n}k^2_n\left(\frac{x-y}{h_n}\right)f(y)\,dy.
\end{gather*}
Then by analogy with (\ref{snb0})
\begin{gather*}
\left(\E_f\|S_n^B\|_2\right)^r=\left\{\E_f\left(\int_{|x|>B}|f_n(x)-\E_f
f_n(x)|^2
dx\right)^{1/2}\right\}^r\\
\leq\frac{L_2^{r/2}}{(n h_n)^{r/2}}\left\{ \|k_n\|_2^2
\int_{|x|>B}f(x)\,dx+ \int_{\Rb}f^{p/2}(x)\,dx
\int_{|x|>B/2h_n}k_n^{2}\left(
x\right)dx\right\}^{r/2}. \\
\end{gather*}
From this, using (\ref{ner11})--(\ref{nerav33}),

\begin{eqnarray*}
\E_f\left(\|S_n^B\|_2^r\right)&\leq& \frac{D^r L_2^{r/2}2^{r-1} r^r}{(1+\log
r)^r}\frac{1}{(n
h_n)^{r/2}}\left[\|k_n\|_2^r\left(\int_{|x|>B}f(x)\,dx\right)^{r/2}
+c_\g h_n\int_{\Rb}f(x)\,dx\right] \\  & &+\frac{D^r L_2^{r/2}2^{2r} r^r h_n^{r/2}\|k_n\|_2^r}{(1+\log
r)^r(n
h_n)^r}
\end{eqnarray*}

\medskip
\n Note that $h_n^{r/2}/(n h_n)^r\to 0$ and, by (\ref{k3}),
$\|k_n\|_2^2\leq\pi^{-1}.$ Therefore for all
$r\geq 1$ there exists a number $n^*>0$
such that for all $n\geq n^*$ and all $f\in\mathbb{F}(\g,2,M)$,

\begin{gather}\label{HHH0}
\E_f\left(\int_{|x|>B}|f_n(x)-\E_f f_n(x)|^2\, dx\right)^{r/2} \leq
\frac{C_{r}}{(n
h_n)^{r/2}}\left(\int_{|x|>B/2}f(x)\,dx\right)^{r/2},
\end{gather}
where $C_{r}$ is a positive constant that may depend on $r$. Because of
the relation (\ref{plus})
the bound (\ref{HHH0}) gives  (\ref{predel1}) when $p=2.$

\bigskip

\n {\textbf{Case 2.}} $1\leq p<2.$ \n We know that (see (\ref{ner22}))

\begin{eqnarray*}
\E_f\|S_n^B\|_p^r &=&\E_f\left(\int_{|x|>B}|f_n(x)-\E_f f_n(x)|^p\, dx\right)^{r/p}\\ &\leq &\frac{D^r 2^{r-1}r^r}{(1+\log r)^r}
\left(\left( \E_f\Vert S^B_{n}\|_p\right) ^{r}+\E_f\max_{1\leq i\leq
n}\Vert Z^B_{i}\Vert_p^{r}\right),
\end{eqnarray*}
\medskip

\n where by (\ref{BBB0}) and (\ref{m3})
\begin{gather*}
\left(\E_f\|S^B_n\|_p\right)^r\leq
\left(\E_f\|S^B_n\|^p_p\right)^{r/p}\leq (n
h_n)^{-r/2}\left(\int_{|x|>B}\left( K_{h_{n}}^{2}\ast f\left(
y\right) \right)^{p/2}dy\right) ^{r/p},
\end{gather*}
and using (\ref{nerav33}) \vspace{-0.5cm}
\begin{gather*}
\E_f\max_{1\leq i\leq n}\|Z_i^B\|^r\leq  \frac{2^{r+1}
h_n^{r/p}\|k_n\|_p^r}{(n h_n)^r}.
\end{gather*}

\medskip

\n Thus, for any $1\leq p<2$ there exists a constant $L_p>0$ such
that for all $r\geq 1,$

\begin{equation}
\E_f\left(\|S_n^B\|_{p}^{r}\right)\leq \left(\frac{r L_{p}}{r+\log
r}\right)^r\left[ \frac{\left( \int_{|x|>B}\left( K_{h_{n}}^{2}\ast
f\left( x\right) \right) ^{p/2}dx\right) ^{r/p}}{\left(
nh_{n}\right) ^{r/2}}+\frac{h_{n}^{r/p}\left\Vert k_n\right\Vert
_{p}^{r}}{\left(
nh_{n}\right)  ^{r}}\right], \label{IB8}%
\end{equation}

\medskip

\n where
$$K^2_{h_n}(x)=h_n^{-1}k_n^2(h_n^{-1}x),\quad x\in\Rb.$$

\bigskip \n Thanks to Lemma 1 of Mason (2009), for all
$n\geq 1$, the integral on the right-hand side of (\ref{IB8}) is bounded
by

\begin{gather}\label{neravenstvo0}
\int_{|x|>B}\left( K_{h_{n}}^{2}\ast f\left( x\right) \right)
^{p/2}dx\leq C{(\lambda,p)} \|k_n\|_2^p
\left\{\E_f(1+|X_{h_n}|^{\lambda})\mathbb{I}(|X_{h_n}|>B)\right\}^{p/2}.
\end{gather}
Here \vspace{-0.5cm}
\begin{gather*}
X_{h_n} \;\;{\bd d\over =}\;\;X+Y_{h_n},
\end{gather*}

\medskip
\n where $X$ has density $f(x)$, $Y_{h_n}$ has density
$K^2_{h_n}(x)/\|k_n\|_2^2$, and $Y_{h_n}$ is independent of $ X$.
Inequality (\ref{neravenstvo0}) implies that

\begin{eqnarray}\label{khn0}
\limsup_{n\to \iy} \int_{|x|>B} \left( K_{h_{n}}^{2}\ast f\left(
x\right) \right) ^{p/2}dx&\leq& C{(\lambda,p)} \|k_n\|_2^p
\left\{\E_f(1+|X|^{\lambda})\mathbb{I}(|X|>B)\right\}^{p/2}\nonumber \\
&\leq&
C_1(\lambda,p)\left(\int_{|x|>B}(1+|x|^{\lambda})f(x)\,dx\right)^{p/2}.
\end{eqnarray}

\medskip
\n The proof of the
upper line in (\ref{khn0}) repeats that of relation (20) in
Lemma 1 of Mason (2009) given for fixed kernel independent of $n$.
We only need to show that for any $(2-p)/p<\l\leq 2$, cf. bound (24) in Mason (2009),

$$\E|Y_{h_n}|^{\l}\to 0,\quad n\to \iy,$$

\n which is true in view of (\ref{eqyn}). Next for all
$1\leq p<2$ and all $r\geq 1$,
\begin{gather*}\frac{h_n^{r/p}\|k_n\|_p^r}{(n h_n)^r}=o\left((n h_n)^{-r/2}\right),\quad n\to \iy.
\end{gather*}

Therefore from the bounds (\ref{IB8}) and  (\ref{khn0}) we get that for any $r\geq 1$
there exists a number $n^*>0$ such that for
all $n\geq n^*$

\begin{gather*}
\E_f\left(\int_{|x|>B}|f_n(x)-\E_f f_n(x)|^p\, dx\right)^{r/p} \leq
\frac{D_{p,r}}{(n
h_n)^{r/2}}\left(\int_{|x|>B}(1+|x|^{\lambda})f(x)\,dx\right)^{r/p},
\end{gather*}

\medskip
\n with some constant $D_{p,r}>0$. In view of condition {\bf (A1)}, this
yields by (\ref{plus}) the statement of Proposition 2 in the case $1\leq p<2.$  The proof
is completed. \done

\bigskip

\n \textbf{Proof of Proposition 3.} Note that by (\ref{ef0}) and (\ref{varf0}), for all $x\in\Rb$, uniformly in $f\in\mathbb{F}(\g,p,M)$,
\begin{gather}\label{disp}
\Var_f\xi_n(x)=\frac{{f(x)}}{\pi}+O(N^{-1}\log N).
\end{gather}
Now let $D_n=D_n(\delta)$ be the region defined as in (\ref{regD}).
When $(x,y)\in D_n$ the random functions $\xi_n(x)$ and $\xi_n(y)$ are weakly correlated and for all large enough $n$,
uniformly in $f\in\mathbb{F}(\g,p,M),$
\begin{eqnarray}\label{xxx}
\sup_{(x,y)\in D_n}\Cov_f(\xi_n(x),\xi_n(y))
=O( N^{-\delta}).
\end{eqnarray}

Indeed, keeping in mind the properties of $f\in\mathbb{F}(\g,p,M)$ (see Section 2.4),
using (\ref{plpl}) and (\ref{ef0}) we obtain that for any $y\in\Rb$, uniformly over $\mathbb{F}(\g,p,M)$,

\vspace{-0.3cm}
\begin{gather*}
\sup_{x\in\Rb}|\Cov_f(k_{h_n}(x-X_1), k_{h_n}(y-X_1))-k_{h_n}(y-x)f(x)|
\end{gather*}
\vspace{-0.6cm}
\begin{gather*} \leq \|k_{h_n}\ast \left[ k_{h_n}(y-\cdot)f(\cdot)\right]-k_{h_n}(y-\cdot)f(\cdot)\|_{\iy}+\|\E^2_f k_{h_n}(\cdot-X_1)\|_{\iy}
\end{gather*}
\vspace{-0.6cm}
\begin{gather*}
\leq
\frac{8M_{\g}n^{1/2}(1+\|k_{h_n}\|_1)}{\pi}\,e^{-\g N}+O(1)=O(\log N).
\end{gather*}

\medskip

\n Thus, for all $x,y\in\Rb$,  uniformly over $ \mathbb{F}(\g,p,M)$,

 \vspace{-0.5cm}
\begin{gather}\label{yyy11}
\Cov_f(k_{h_n}(x-X_1), k_{h_n}(y-X_1))=k_{h_n}(y-x)f(x)+O(\log N),
\end{gather}

\medskip
\n where for some constant $C=C({\g,p,M})>0$,

\begin{gather*}|k_{h_n}(y-x)f(x)|
 =\frac{N h_n
f(x)}{\pi(x-y)^2} \left|\cos(N(x-y))-\cos(h_n^{-1}(x-y))\right|\leq
\frac{C}{(x-y)^{2}}.\end{gather*}
\medskip

\n The last inequality implies that  for all large enough $n$, uniformly over ${\mathbb{ F}}(\g,p,M)$,

$$\sup_{(x,y)\in D_n}|k_{h_n}(y-x)f(x)|=O( N^{1-\delta}).$$

\n Therefore, noting that \begin{gather*}
\Cov_f(\xi_n(x),\xi_n(y))=h_n\Cov_f(k_{h_n}(x-X_1), k_{h_n}(y-X_1)),
\end{gather*}
where $h_n=N^{-1}(1-N^{-1})$, and using (\ref{yyy11}) we get (\ref{xxx}).

 We next turn to the sequence $(\xi_n(x), \xi_n(y))=\sum_{i=1}^n \etab_{i,n}(x,y)$, $n=1,2,\ldots,$ where $\etab_{i,n}(x,y)=(\eta_{i,n}(x), \eta_{i,n}(y))$ is given by (\ref{etab}), and show
that the Lindeberg condition
\begin{gather}\label{Lc}
\lim_{n\to \iy}\sup_{(x,y)\in D_n}n\E_f\left(|\etab_{1,n}|^2\mathbb{I}(|\etab_{1,n}|>\tau)\right)=
0,\quad \forall\; \tau>0,
\end{gather}
holds uniformly in $f\in\mathbb{F}(\g,p,M)$.
To this end observe that by (\ref{varf0}) and Chebyshev's
 inequality,  uniformly in $x\in\Rb$ and $f\in\mathbb{F}(\g,p,M)$,
 \begin{gather*}
\Pb_f\left((\eta_{1,n}(x))^2>{\tau^2}/{2}\right)=
O(n^{-1}).
\end{gather*}
 Therefore by (\ref{fmax0}) and (\ref{ef0})  for all sufficiently large $n$ and some constant $C=C({\g,p,M})>0$,  uniformly over $\mathbb{F}(\g,p,M)$,
\begin{gather*}
 \sup_{(x,y)\in D_n}
{n}\E_f\left(|\etab_{1,n}|^2 \mathbb{I}(|\etab_{1,n}|>\tau)\right)\\
\leq {C}(n h_n)^{-1}\sup_{(x,y)\in D_n}n \Pb_f\left((\eta_{n,1}(x))^2+(\eta_{n,1}(y))^2>\tau^2\right)\\
  \leq
{2C}{(nh_n)^{-1}} \sup_{x\in\Rb} n \Pb_f\left((\eta_{1,n}(x))^2>{\tau^2}/{2}\right)
\\ =
O\left((n h_n)^{-1}\right)=o(1).
\end{gather*}
This ensures the validity of (\ref{Lc}).

Finally, using Fact 1 we infer from (\ref{disp}), (\ref{xxx}), and (\ref{Lc}) that the statement of Proposition 3 holds. \done

\section*{Acknowledgements} This research was partly supported by an NSERC
grant. We would like to thank Boris Levit for helpful discussions and suggestions.

\end{document}